\documentclass[reqno]{amsart}
\usepackage{amssymb}
\usepackage{amsmath}

\makeatletter
\@addtoreset{equation}{section}
\makeatother

\renewcommand\thefigure{\thesection.\@arabic\c@figure}
\renewcommand\thetable{\thesection.\@arabic\c@table}

\newtheorem{theorem}{Theorem}[section]
\newtheorem{lemma}[theorem]{Lemma}
\newtheorem{proposition}[theorem]{Proposition}
\newtheorem{corollary}[theorem]{Corollary}

\newcommand{\mc}[1]{{\mathcal #1}}
\newcommand{\mf}[1]{{\mathfrak #1}}
\newcommand{\mb}[1]{{\mathbf #1}}
\newcommand{\bb}[1]{{\mathbb #1}}

\def\emptysquare{{\hbox{\vrule height6pt width0.6pt depth0pt%
\vbox{\hrule height0.6pt width4.8pt depth0pt%
\vglue4.8pt%
\hrule height0.6pt width4.8pt depth0pt}%
\vrule height6pt width0.6pt depth0pt}}}

\def\qed{\unskip\nobreak
\hfil\penalty50\hskip1.75em\null\nobreak\hfil\emptysquare
{\parfillskip=0pt \finalhyphendemerits=0 \par}\medskip}

\begin{document}

\title[Symmetric exclusion process]{Gaussian estimates for symmetric
  simple exclusion processes.}

\author{C. Landim}   

\address{\noindent IMPA, Estrada Dona Castorina 110,
CEP 22460 Rio de Janeiro, Brasil and CNRS UMR 6085,
Universit\'e de Rouen, 76128 Mont Saint Aignan, France. 
\newline
e-mail:  \rm \texttt{landim@impa.br}
}

\subjclass[2000]{primary 60K35; secondary 82A05} 

\keywords{Interacting particle systems, decay to equilibrium,
  hydrodynamic equation, fluctuations of density field}

\thanks{}

\dedicatory{In memoriam of Martine Babillot.}

\begin{abstract}
  We prove Gaussian tail estimates for the transition probability of
  $n$ particles evolving as symmetric exclusion processes on $\bb
  Z^d$, improving results obtained in \cite{l}. We derive from this
  result a non-equilibrium Boltzmann-Gibbs principle for the symmetric
  simple exclusion process in dimension $1$ starting from a product
  measure with slowly varying parameter.
\end{abstract}

\maketitle                

\section{Introduction}
\label{sec0}

To derive sharp bounds on the rate of convergence to equilibrium is
one of the main questions in the theory of Markov processes. In the
last decade, this problem has attracted many attention in the context
of conservative interacting particle systems in infinite volume. Fine
estimates of the spectral gap of reversible generators restricted to
finite cubes and logarithmic Sobolev inequalities have been obtained.
We refer to \cite{l} for the recent literature on the subject. From
these bounds on the ergodic constants, polynomial decay to equilibrium
in $L^2$ has been proved for some processes.  For instance, Bertini
and Zegarlinski \cite{bz1}, \cite{bz2} proved that the symmetric
simple exclusion process in $\bb Z^d$ converges to equilibrium in
$L^2$ at rate $t^{-d/2}$. For a class of functions $f$ that includes
the cylinder functions, there exists $V(f)$ finite such that
$$
\Vert P_t f - <f>_\alpha \Vert_2^2 \; \le\; \frac{V(f)}{(1+t)^{d/2}}
$$
for all $t\ge 0$. Here $P_t$ stands for the semi--group, $<f>_\alpha$
for the expectation of $f$ with respect to the Bernoulli product
measure with density $\alpha$ and $\Vert f\Vert_2$ for the $L^2$ norm of
$f$.  Janvresse, Landim, Quastel and Yau \cite{jlqy} and Landim and
Yau \cite{ly} extended the algebraic decay in $L^2$ for zero range
and Ginzburg-Landau dynamics.
                                                                               
We refine in this article the Gaussian upper bounds obtained in
\cite{l} for the transition probabilities of finite symmetric simple
exclusion processes evolving on the lattice $\bb Z^d$. Our approach is
based on a logarithmic Sobolev inequality and on Davies \cite{d}
method to derive estimates for heat kernels.

Consider $n\ge 2$ indistinguishable particles moving on the
$d$-dimensional lattice $\bb Z^d$ as symmetric random walks with an
exclusion rule which prevents more than one particle per site. The
dynamics can be informally described as follows. There are initially
$n$ particles on $n$ distinct sites of $\bb Z^d$. Each particle waits
a mean one exponential time at the end of which, being at $x$, it
chooses a site $y$ with probability $p(y-x)$, for some finite range,
irreducible, symmetric transition probability $p(\cdot)$. If the site
is vacant, the particle jumps, otherwise it stays where it is and
waits a new mean one exponential time.

The state space of this Markov process, denoted by $\mc E_{n}$, is the
collection of all subsets $A$ of $\bb Z^d$ with cardinality $n$:
$$
\mc E_{n} \;=\; \{A \subset \bb Z^d\,, |A| = n\}\;;
$$
while its generator $\mc L_n$ is given by
\begin{equation}
\label{eq:9}
(\mc L_nf)(A) \;=\; \sum_{x,y\in\bb Z^d} p(y-x)
[f(A_{x,y}) - f(A)] \; ,
\end{equation}
where $A_{x,y}$ stands for the set $A$ with sites $x$, $y$ exchanged:
$$
A_{x,y} \;=\;
\left\{
\begin{array}{cl}
(A \setminus \{x\}) \cup \{y\} & \text{if $ x \in A$, $y \not\in
  A$,} \\
(A \setminus \{y\}) \cup \{x\} & \text{if $ y \in A$, $x \not\in
  A$,} \\
A  & \text{otherwise.}
\end{array}
\right.
$$
In formula (\ref{eq:9}) summation is carried over all bonds $\{y,z\}$
to avoid counting twice the contribution of the same jump.

It is easy to check that the counting measure on $\mc E_{n}$\,,
denoted by $\nu_{n}$\,\, ($\nu_{n}(A) = 1$ for every $A$ in $\mc
E_{n}$), is an invariant, reversible measure for the process.

Fix a set $A_0$ in $\mc E_n$ and denote by $f_t = f_t^{A_0}$ the
solution of the forward equation with initial data $\delta_{A_0}$~:
\begin{equation}
\label{eq:10a}
\left\{
\aligned
& \partial_t f_t = \mc L_n f_t \; , \\
& f_0 (A) = \mb 1\{ A = A_0\} \; .
\endaligned
\right.  
\end{equation}

The main result of the article provides a Gaussian estimate for the
transition probability $f_t$. Denote by $\mb x = (x_1, \dots, x_n)$
the sites of $(\bb Z^d)^n$. For a configuration $\mb x$, let $\mb
x_{i,j}$ be the $j$-th coordinate of the $i$-th point of $\mb x$: $\mb
x_{i,j} = x_i\cdot e_j$, where $\cdot$ stands for the inner product in
$\bb R^d$ and $\{e_1, \dots, e_d\}$ for the canonical basis of $\bb
R^d$. The Euclidean norm of $(\bb R^d)^n$ is denoted by $\Vert \mb
x\Vert$ so that $\Vert \mb x\Vert^2 = \sum_{i,j} x_{i,j}^2$.

Denote by $\Phi$ the Legendre transform of the convex function $w^2
\cosh w$:
$$
\Phi(u) \;=\; \sup_{w\in\bb R} \big\{ u w - w^2
\cosh w \big\}\;.
$$
An elementary computations shows that $\Phi (w) \sim w^2$ for $w$
small and $\Phi (w) \sim w \log w$ for $w$ large. 

\begin{theorem}
\label{s5}
Fix a set $A_0 = \{z_1,\dots,z_n\}$ in $\mc E_n$\,.  Let $f_t$ be the
solution of the forward equation (\ref{eq:10a}). There exist finite
constants $C_2 = C_2(n,d,p)$, $a_0=a_0(p)$ such that
$$
f_T(A) \;\le\, \sum_\sigma \frac{C_2}{(1+T)^{nd/2}}\,
\exp\Big\{- \frac { a_0 T}{2 (\log T)^2} \, \Phi\Big( 
\frac{\Vert \mb x_\sigma -\mb z \Vert \log T} {a_0^2  T}\Big)\Big\}
$$
for every $T > C_2$ and every set $A = \{x_1,\dots,x_n\}$. In this
formula, summations is performed over all permutations $\sigma$ of $n$
points and $\mb x_\sigma$ stands for the vector $(x_{\sigma (1)}, \dots
  ,x_{\sigma (n)})$.
\end{theorem}

The asymptotic behavior of $\Phi(\cdot)$ at the origin shows that for
every $\gamma >0$, there exists a constant $a_1 = a_1(p,\gamma)$ such
that
\begin{equation}
\label{eq:10}
f_T(A) \;\le\, \sum_\sigma \frac{C_2}{(1+T)^{nd/2}}\,
\exp\Big\{- \frac { \Vert \mb x_\sigma -\mb z \Vert^2}{a_1 T} \Big\}
\end{equation}
for all $T>C_2$ and all sets $A$ such that $\Vert \mb x_\sigma -\mb z
\Vert \le \gamma T/\log T$ for all permutations $\sigma$. Furthermore,
since
$$
\Phi\Big( \frac{\Vert \mb x_\sigma -\mb z \Vert \log T} {a_0^2
  T}\Big) \;\ge\; \frac 1n \sum_{i=1}^n \Phi\Big( \frac{\Vert \mb
  x_{\sigma(i)} - \mb z_i \Vert \log T} {a_0^2 T}\Big)
$$
we have that
$$
f_T(A) \;\le\, \sum_\sigma \frac{C_2}{(1+T)^{nd/2}}\,
\exp\Big\{- \frac { a_0 T}{2n (\log T)^2} \, \sum_{i=1}^n \Phi\Big( 
\frac{\Vert \mb x_{\sigma(i)} -\mb z_i \Vert \log T} {a_0^2
  T}\Big)\Big\}\; .
$$
For a fixed $\gamma>0$, in last formula we may replace $\Phi(w)$ by
$C(\gamma)w^2$ if $\Vert \mb x_{\sigma(i)} -\mb z_i \Vert \le a_0^2
T \gamma/\log T$ and $\Phi(w)$ by $C(\gamma)w \log w $ otherwise.

\section{Boltzmann-Gibbs principle}
\label{sec2}

We prove in this section the Boltzmann-Gibbs principle for the
symmetric simple exclusion process out of equilibrium in dimension 1.
This result allows the replacement of average of local functions by
functions of the empirical density in the fluctuation regime and is
the main point in the proof of a central limit theorem around the
hydrodynamical limit for interacting particle systems (cf. \cite{kl},
Chap. 11). We restricted our attention to dimension 1 because Lemma
\ref{s10} below has only been proved in $d=1$. 

The Boltzmann-Gibbs principle for one-dimensional processes out of
equilibrium was proved in \cite{cy} through the logarithmic Sobolev
inequality. A different version, involving microscopic time integrals,
is presented in \cite{fpv} and uses sharp estimates on the comparison
between independent random walks and the symmetric exclusion process.

Fix a profile $\rho_0: \bb R^d \to [0,1]$ in $C^1(\bb R^d)$ with a
bounded derivative and consider a sequence of product measures
$\{\nu^N_{\rho_0(\cdot)},\, N\ge 1\}$ on $\{0,1\}^{\bb Z^d}$
associated to this profile so that
$$
E_{\nu^N_{\rho_0(\cdot)}}[\eta(x)]\;=\; \rho_0(x/N)\;.
$$
Denote by $\bb P_{\nu^N_{\rho_0(\cdot)}}$ the probability measure
on the path space $D(\bb R_+, \{0,1\}^{\bb Z^d})$ corresponding to the
symmetric simple exclusion process starting from
$\nu^N_{\rho_0(\cdot)}$ and speeded up by $N^2$.  Expectation with
respect to $\bb P_{\nu^N_{\rho_0(\cdot)}}$ is denoted by $\bb
E_{\nu^N_{\rho_0(\cdot)}}$.

For $x$ in $\bb Z^d$, let
$$
\rho^N(t,x) \;=\; \bb E_{\nu^N_{\rho_0(\cdot)}} [ \eta_t(x)]\; .
$$
Of course, $\rho^N(t,x)$ is the solution of the linear equation
$$
\left\{
\aligned
& \partial_t \rho^N(t,x) = N^2 \sum_{y\in\bb Z^d} p(y-x) 
[\rho^N(t,y) - \rho^N(t,x)]  \; , \\
&  \rho^N(0,x) = \rho_0(x/N) \;.
\endaligned
\right.  
$$
This equation can be written as $\partial_t \rho^N(t,x) = N^2 \mc
L_1 \rho^N(t,x)$, where $\mc L_1$ is the generator introduced in
(\ref{eq:9}). Next proposition is the main result of this section.

\begin{proposition}
\label{s7}
Let $d=1$ and fix $T>0$, a finite subset $A$ of $\bb Z$ such that
$|A|>2$ and a continuous function $H$ in $L^1(\bb R)$. Then,
$$
\lim_{N\to\infty} \bb E_{\nu^N_{\rho_0(\cdot)}} \Big[\Big( \int_0^T
dt \, \frac 1{N^{1/2}} \sum_{x\in\bb Z} H(x/N) \prod_{z\in A}
[\eta_t (x+z) - \rho^N(t,x+z)]\Big)^2 \Big]\;=\;0\;.
$$
\end{proposition}

The proof of this proposition is presented at the end of this
section. The Boltzmann-Gibbs principle is a simple
consequence but requires some extra notation.

For a finite subset $A$ of $\bb Z$ and $0\le \alpha\le 1$, let 
$$
\Psi (A,\alpha) \;=\; \prod_{x\in A} [\eta (x) - \alpha]\;.
$$
By convention, we set $\Psi (\phi,\alpha) = 1$.  Each cylinder
function $f:\{0,1\}^{\bb Z} \to \bb R$ can be written as
$$
f(\eta) \;=\; \sum_{A\in\mc E} \mf f(A, \alpha) \Psi (A,\alpha)\;.
$$
A straightforward computation shows that for each finite set $A$,
$\mf f(A, \cdot)$ is a smooth function, in fact a polynomial.

For a cylinder function $f$, let $\tilde f :[0,1]\to\bb R$ be the real
function defined by $\tilde f (\alpha) = E_{\nu_\alpha} [f(\eta)]$ and
let
$$
\Gamma_f(\eta,\alpha)\;=\; f(\eta) \;-\; \tilde f (\alpha) \;-\;
\tilde f'(\alpha) [\eta (0) - \alpha]\;.
$$
Note that $\mf f(\phi, \alpha)= \tilde f(\alpha)$ and that
$\sum_{x\in\bb Z} \mf f(\{x\}, \alpha) = \tilde f'(\alpha)$. In
particular, it follows from a simple computation that
$$
\Gamma_f(\eta,\alpha)\;=\; \sum_{z\in \bb Z} \mf f(\{z\}, \alpha) [
\eta (z) - \eta(0)] \;+\; \sum_{n\ge 2}\sum_{A\in \mc E_n} \mf f(A,
\alpha) \Psi (A,\alpha)\;.
$$

Fix a smooth functions $H$ in $L^1(\bb R)$. By the previous formula,
\begin{eqnarray*}
\!\!\!\!\!\!\!\!\!\!\!\!\!\! &&
\frac 1{N^{1/2}} \sum_{x\in\bb Z} H(x/N) 
\Gamma_f(\tau_x \eta, \rho^N(t,x)) \\
\!\!\!\!\!\!\!\!\!\!\!\!\!\! && \quad \;=\;
\sum_{z\in \bb Z} \frac 1{N^{1/2}}  \sum_{x\in\bb Z} H(x/N) 
\mf f(\{z\}, \rho^N(t,x))  [\eta (x+z) - \eta(x)] \\
\!\!\!\!\!\!\!\!\!\!\!\!\!\! && \qquad
+\; \sum_{n\ge 2}\sum_{A\in \mc E_n} \frac 1{N^{1/2}} \sum_{x\in\bb Z}
H(x/N) \mf f(A, \rho^N(t,x)) \tau_x \Psi (A, \rho^N(t,x))\;.
\end{eqnarray*}
Note that the sums in $z$, $n$ and $A$ are finite because $f$ is a
cylinder function. Since $H$, $\rho^N(t,\cdot)$ and $\mf f(\{z\}, \cdot)$
are smooth functions, a change of variables shows that the first term
is of order $N^{-1/2}$ and that the second is equal to
$$
\sum_{n\ge 2}\sum_{A\in \mc E_n} \frac 1{N^{1/2}} \sum_{x\in\bb Z}
H(x/N) \mf f(A, \rho^N(t,x)) \prod_{z\in A} [\eta (x+z) - \rho^N(t,x+z)]
\;+\; O(N^{-1/2})\; ,
$$
which is exactly the expression appearing in Proposition \ref{s7}.
Since $\mf f(A, \cdot)$ and $\rho^N(t,\cdot)$ are smooth bounded
functions, we have proved the following result, known as the
Boltzmann-Gibbs principle.

\begin{corollary}
\label{s11}
Let $d=1$ and fix $T>0$, a cylinder function $f$ and a smooth function
$H$ in $L^1(\bb R)$. Then,
$$
\lim_{N\to\infty} \bb E_{\nu^N_{\rho_0(\cdot)}} \Big[\Big( \int_0^T
dt \, \frac 1{N^{1/2}} \sum_{x\in\bb Z} H(x/N) 
\Gamma_f(\tau_x \eta_t, \rho^N(t,x)) \Big)^2 \Big]\;=\;0\;.
$$
\end{corollary}

The proof of Proposition \ref{s7} is based on three lemmas concerning
the decay of the space-time correlations of the symmetric exclusion
process.  We start with a general result which will be used
repeatedly. 

Fix $n\ge 2$ and denote by $f_t(A,B)=f^N_t(A,B)$ the semi-group
associated to the generator $N^2 \mc L_n$.  For a finite subset $A$ of
$\bb Z$, let
$$
I(A) \;=\; \sum_{x,y\in A} p(y-x)\; . 
$$
Note that $I(A)$ vanishes unless $A$ contains two sites which are
within a distance smaller than the range of the transition
probability. Next lemma follows from Theorem \ref{s5} and a
straightforward computation.

\begin{lemma}
\label{s6}
For all $T<\infty$, $n\ge 2$, there exists a finite constant $C_3$,
depending only on $n$, $p$ and $T$ such that
$$
\int_0^t ds\; \frac 1{(1+sN^2)^{m/2}}\sum_{B\in\mc E_n} f^N_{t-s}(A, B) \,
I(B) \;\le\; C_3 A_N (m,t) 
$$
for all $A$ in $\mc E_n$, $N\ge 1$, $0\le t\le T$; where
$$
A_N(m,t) = 
\left\{
  \begin{array}{cl}
N^{-1} & \text{if $m=0$,} \\
N^{-2} & \text{if $m=1$,} \\
\log N/ N^2 \sqrt{1+tN^2} & \text{if $m=2$,} \\
1/ N^2 \sqrt{1+tN^2} & \text{if $m\ge 3$.}
  \end{array}
\right.
$$
\end{lemma}

We now introduce the space-time correlations, also called
$v$-functions in \cite{fpsv}.  For a finite subset $A$ of $\bb Z$
and $t\ge 0$, let $\varphi^N(t,\phi)=1$,
$$
\varphi^N(t,A)\;=\; \bb E_{\nu^N_{\rho_0(\cdot)}} 
\Big[ \prod_{x\in A} \{\eta_t(x) - \rho^N(t,x)\}  \Big]\; .
$$
Notice that $\varphi^N(t,\{x\})$ vanishes for all $x$. An
elementary computation shows that
\begin{equation}
\label{eq:11a}
\left\{
\aligned
&  \partial_t \varphi^N(t,A)\; =\; N^2 \mc L_n \varphi^N(t,A) 
\;+\; G^N(t,A) \; , \\
&  \varphi^N(0,A) = 0 \;,
\endaligned
\right.
\end{equation}
where $n=|A|$ and $G^N(t,A)$ is given by
\begin{eqnarray*}
\!\!\!\!\!\!\!\!\!\!\!\!\!\!\!\! &&
N^2 \sum_{y,z \in A } p(z-y) \Big\{ \varphi^N(t,A\setminus \{z\})
- \varphi^N(t,A\setminus \{y\}) \Big\} 
\Big\{ \rho^N(t,\{z\}) - \rho^N(t,\{y\}) \Big \} \\
\!\!\!\!\!\!\!\!\!\!\!\!\!\!\!\! && \quad
-\; (N^2/2) \sum_{y,z \in A } p(z-y) \varphi^N(t,A\setminus \{y,z\}) 
\Big\{ \rho^N(t,\{z\}) - \rho^N(t,\{y\}) \Big \}^2 \;.
\end{eqnarray*}
Here again summation is carried over all bonds.  Notice that the first
line vanishes for $n=2$ and that the second line vanishes for
$n=3$.

The linear differential equation (\ref{eq:11a}) has a unique solution
which can be represented as
$$ 
\varphi^N(t,A)\;=\; \int_0^t ds\, \sum_{B\in\mc E_n} f_{t-s} (A,B)
G^N(s,B)
$$
so that the space-time correlations $\varphi^N(t,A)$ can be
estimated inductively in $n$. 

Next lemma is due to Ferrari, Presutti, Scacciatelli and Vares
\cite{fpsv}. In the proof of Proposition \ref{s7} we don't need such
sharp estimates. 

\begin{lemma}
\label{s10}
Assume that $d=1$ and fix $T>0$. For each $n\ge 1$, there exists a
finite constant $C_4=C_4(n,p,\rho_0,T)$ such that
\begin{equation*}
\sup_{\substack{0\le t\le T\\ A\in\mc E_{2n}}}
\big\vert \varphi^N(t,A) \big\vert \;\le\; \frac{C_4}{N^n} \;, \quad
\sup_{\substack{0\le t\le T\\ A\in\mc E_{2n+1}}}
\big\vert \varphi^N(t,A) \big\vert \;\le\; \frac{C_4 \log N }{N^{n+1}}\;.
\end{equation*}
\end{lemma}
\medskip

For $0\le s\le t\le T$, $1\le k\le n$ and $B\in \mc E_k$, $A\in\mc
E_n$, let
$$
R_N(s,A;t,B) \;=\; \bb E_{\nu^N_{\rho_0(\cdot)}} 
\Big[\prod_{x\in A} [\eta_s (x) - \rho^N(s,x)] 
\prod_{y\in B} [\eta_t (y) - \rho^N(t,y)] \, \Big]\;.
$$
Since $s$ and $A$ will be fixed, most of the time, we denote
$R_N(s,A;t,B)$ by $R_N(t,B)$. Notice that in this definition we do not
require $A$ and $B$ to have the same cardinality.  An elementary
computation shows that $R_N(t,B)$ is the solution of the linear
differential equation
\begin{equation}
\label{eq:11}
\left\{
\aligned
&  \partial_t R_N(t,B)\; =\; N^2 \mc L_k R_N(t,B) \;+\; H_N(t,B)
\quad \text{for $t\ge s$}\; , \\
&  R_N(s,B) = J_N(s,A,B) \;,
\endaligned
\right.
\end{equation}
where $k=|B|$, 
$$
J_N(s,A,B)\;=\; \bb E_{\nu^N_{\rho_0(\cdot)}} 
\Big[\prod_{x\in A} [\eta_s (x) - \rho^N(s,x)] 
\prod_{y\in B} [\eta_s (y) - \rho^N(s,y)] \, \Big]
$$
and $H_N(t,B)$ is given by
\begin{eqnarray*}
\!\!\!\!\!\!\!\!\!\!\!\!\!\! &&
N^2 \sum_{y,z \in B } p(z-y) \Big\{ R_N(t,B\setminus \{z\})
- R_N(t,B\setminus \{y\}) \Big\} 
\Big\{ \rho^N(t,\{z\}) - \rho^N(t,\{y\}) \Big \} \\
\!\!\!\!\!\!\!\!\!\!\!\!\!\! && \quad 
-\; (N^2/2) \sum_{y,z \in B } p(z-y) R_N(t,B\setminus \{y,z\})
\Big\{ \rho^N(t,\{z\}) - \rho^N(t,\{y\}) \Big \}^2\; .
\end{eqnarray*}
Notice that $R_N(t,\phi) = \varphi^N(s,A)$, that $H_N(t,B)$ vanishes
for $n=1$ and that $J_N(s,A,B)$ is not equal to $\varphi^N(s, A\cup
B)$ but given by
$$
\sum_{C} \varphi^N(s, (A \Delta B)\cup
C) \prod_{x\in C} [1-2 \rho^N(s,\{x\})] \prod_{x\in (A \cap B)
  \setminus C} \rho^N(s,\{x\}) [1-\rho^N(s,\{x\})]\;,
$$
where the summation is carried over all subsets $C$ of $A\cap B$ and
where $A \Delta B$ stands for the symmetric difference of $A$ and $B$.

The differential equation (\ref{eq:11}) has a unique solution which
can be represented as
\begin{equation}
\label{eq:14}
R_N(t,B)\;=\; \sum_{C\in\mc E_k} f_{r} (B,C)\, J_N(s,A,C)
\;+\; \int_0^{r} du\, \sum_{C\in\mc E_k} f_{r-u} (B,C)
H_N(s+u,C)\; ,
\end{equation}
where $r=t-s$. This last notation is systematically used below. Let
$$
U_N(t,B)\;=\; \sum_{C\in\mc E_k} f_{r} (B,C)\, J_N(s,A,C)\; .
$$

\begin{lemma}
\label{s12}
Fix $2\le k\le n$, $0\le s\le t\le T$ and $A$ in $\mc E_n$. There
exists a finite constant $C_4=C_4(p,n,T,\rho_0)$ such that
$$
\sup_{B\in\mc E_k} \big\vert U_N(t, B)\big\vert \;\le\; 
\frac{C_4 B(n-k)}{(1+rN^2)^{k/2}}\;,
$$
where $B(2j) = N^{-j}$ and $B(2j+1)= \log N / N^{j+1}$ for $j\ge
0$.
\end{lemma}

\begin{proof}
$U_N(t,B)$ is absolutely bounded by
\begin{equation}
\label{eq:13}
\sum_{j=0}^k J(s,j) \sum_{\substack{D\subset A, D\in\mc E_j\\
E\cap A=\phi, E\in \mc E_{k-j}}} f_r(B, D\cup E)\;,
\end{equation}
where
$$
J(s,\ell)\;=\; \sup_{\substack{C\in \mc E_k\\ |C\cap A|=\ell}} 
\big\vert J_N(s, A, C) \big\vert \;\le\; \max_{m}
\sup_{D\in \mc E_m} \big\vert \varphi^N(s, D) \big\vert\; ,
$$
where the maximum is carried over $n+k-2\ell \le m\le n+k+\ell$.
Last inequality follows from the explicit formula for $J_N(s,A,C)$. 
By Lemma \ref{s10}, the previous expression is less than or equal
to $C_4 B(n+k-2\ell)$.  On the other hand, by Theorem \ref{s5},
$$
\sum_{\substack{D\subset A, D\in\mc E_j\\
E\cap A=\phi, E\in \mc E_{k-j}}} f_r(B, D\cup E)\;\le\;
\frac {C_2(k,p)}{(1+rN^2)^{j/2}}\; \cdot
$$
Therefore,
$$
\big\vert U_N(t, B)\big\vert \;\le\; C_4 \sum_{j=0}^k 
\frac{B(n+k-2j)}{(1+rN^2)^{j/2}}\;\le\;  
\frac{C_4 B(n-k)}{(1+rN^2)^{k/2}}\;\cdot 
$$
This concludes the proof of the lemma.
\end{proof}

We are now in a position to prove the main result towards the
Boltzmann-Gibbs principle.

\begin{lemma}
\label{s9}
Fix $n\ge 2$, there exists a finite constant $C_4=C_4(n,p,T,\rho_0)$
such that
\begin{equation*}
\big\vert R_N(s,A;t,B)\big\vert  \;\le\; C_4\Big\{
\frac {\log N}{N^2} + \frac{1}{1+ (t-s)N^2} \Big\}
\end{equation*}
for all $A$, $B$ in $\mc E_n$, $0\le s<t\le T$.
\end{lemma}

\begin{proof}
Fix $n\ge 2$, $s\ge 0$ and $A$ in $\mc E_n$. For $1\le k\le n$, denote
by $R_{N,k}(t, \cdot)$ the solution of the linear differential
equation (\ref{eq:11}). Since the equation for $R_{N,k}$ involves
$R_{N,k-1}$, $R_{N,k-2}$, an induction argument on $k$ is required.  A
simple pattern appears only for $k\ge 7$. Hence, for $1\le k\le 6$, we
need to proceed by inspection, making the proof long and tedious.

Consider $k=1$. In this case $H_N$ vanishes and, by Lemma
\ref{s12},
$$
\big\vert R_{N,1} (t, \{x\}) \big\vert \;=\; 
\big\vert U_{N} (t, \{x\}) \big\vert \;\le\; 
\frac{a_1 B(n-1)}{(1+rN^2)^{1/2}} \;\cdot
$$
Here and below $\{a_j, j\ge 1\}$ are finite constants depending on
$n$, $p$, $T$ and $\rho_0$ which may change from line to line.

For $k=2$, since $R_N(t,\phi)$ is time independent and absolutely
bounded by $B(n)$, the previous estimates and Lemma \ref{s12} show
that  
$$
\big\vert H_{N,2} (t, B) \big\vert \;\le\; 
a_2\Big\{ \frac{N B(n-1)}{(1+rN^2)^{1/2}} \;+\; B(n) \Big\} I(B) \; .
$$
Therefore, by the explicit formula (\ref{eq:14}) for $R_{N,2}(t,B)$
and by Lemmas \ref{s6}, \ref{s12},
$$
\big\vert R_{N,2} (t, B) \big\vert \;\le\; 
a_2\Big\{ \frac{B(n-2)}{1+rN^2} \;+\; \frac{B(n-1)}{N} \Big\}  
$$
because $B(n)\le B(n-1)$. Notice that this inequality proves the lemma
for $n=2$ because $R_{N,2} (t, B) = R_N(s,A;t,B)$.

The estimates for $R_{N,1}$ and $R_{N,2}$ give bounds for $H_{N,3}$
which in turn, together with the explicit formula (\ref{eq:14}) for
$R_{N,3}(t,B)$ and Lemmas \ref{s6}, \ref{s12} show that
$$
\big\vert R_{N,3} (t, B) \big\vert \;\le\; 
a_3\Big\{ \frac{B(n-3)}{(1+rN^2)^{3/2}} \;+\; \frac{B(n-2) \log
N}{N(1+rN^2)^{1/2}} \Big\}   \;.
$$
Here we used the fact that $B(n-3) = N B(n-1)$ to eliminate one of the
terms appearing in the expression of $R_{N,3} (t, B)$.

We repeat this procedure for $k=4$, $5$ and $6$ to obtain that
\begin{eqnarray*}
\!\!\!\!\!\!\!\!\!\!\!\!\!\!\! &&
\big\vert R_{N,4} (t, B) \big\vert \;\le\; 
a_4\Big\{ \frac{B(n-4)}{(1+rN^2)^{2}} \;+\; \frac{B(n-3)}
{N(1+rN^2)^{1/2}} \;+\; \frac{B(n-2) \log N} {N^{2}}  \Big\}   \;, \\
\!\!\!\!\!\!\!\!\!\!\!\!\!\!\! &&
\quad \big\vert R_{N,5} (t, B) \big\vert \;\le\; 
a_5\Big\{ \frac{B(n-5)}{(1+rN^2)^{5/2}} \;+\; \frac{B(n-4)}
{N(1+rN^2)^{1/2}}  \;+\; \frac{B(n-3)} {N^{2}}\Big\}   \;, \\
\!\!\!\!\!\!\!\!\!\!\!\!\!\!\! && \qquad 
\big\vert R_{N,6} (t, B) \big\vert \;\le\; 
a_6\Big\{ \frac{B(n-6)}{(1+rN^2)^{3}} \;+\; \frac{B(n-5)}
{N(1+rN^2)^{1/2}} \;+\; \frac{B(n-4)} {N^{2}}  \Big\}   \;\cdot
\end{eqnarray*}
For $k=5$, we used the fact that $B(\ell) \log N \le B(\ell-1)$.

A pattern has been found for $k=5$, $6$. It is now a simple matter to
prove by induction that this pattern is conserved so that
$$
\big\vert R_{N,k} (t, B) \big\vert \;\le\; 
a_6\Big\{ \frac{B(n-k)}{(1+rN^2)^{k/2}} \;+\; \frac{B(n-k+1)}
{N(1+rN^2)^{1/2}} \;+\; \frac{B(n-k+2)} {N^{2}}  \Big\} 
$$
for $k\ge 7$. It remains to recall the definition of $B(j)$ and to
recollect all previous estimates to conclude the proof of the lemma.
\end{proof}

Notice that we could have set $B(1) = N^{-1}$ for the estimates in the
previous lemma. Taking $B(1) = \log N/N$ simplifies slightly the
notation since we have that $B(n+2) = B(n) N^{-1}$ for all $n\ge 0$
and we miss only a $\log N$ factor, which is irrelevant for our
purposes. 

We are now in a position to prove Proposition \ref{s7}.  With the
notation introduced in this section, the expectation appearing in the
statement of the proposition becomes
$$
\frac 2N \sum_{x,y\in\bb Z} H(x/N)\, H(y/N) \int_0^T dt\, 
\int_0^t ds\, R_N(s,A+x;t,A+y)\;,
$$
where $A+x$ is the set $\{z+x:\, z\in A\}$. By Lemma \ref{s9} and a
change of variables, this expression is bounded above by
\begin{eqnarray*}
\frac {C(n,p,\rho_0,T) \log N}N \Big( \frac 1N \sum_{x\in\bb Z} 
\big\vert H(x/N) \big\vert \Big)^2 \;,
\end{eqnarray*}
which proves Proposition \ref{s7}.

We conclude this section with an observation. The same arguments
presented above in the proof of Proposition \ref{s7} shows that
\begin{eqnarray*}
\!\!\!\!\!\!\!\!\!\!\!\!\!\!\! &&
\bb E_{\nu^N_{\rho_0(\cdot)}} \Big[\Big( \int_0^T
dt \, \frac 1{N^{1/2}} \sum_{x\in\bb Z} H(x/N) 
[\eta_t (x) - \rho^N(t,x)]\Big)^2 \Big] \\
\!\!\!\!\!\!\!\!\!\!\!\!\!\!\! && \quad
\le \; 2 \int_0^T dt\, \int_0^t ds\,  
\frac 1N \sum_{x\in \bb Z} F(\rho^N (s,x))\, H(x/N)\, 
(f^N_{t-s} H)(x/N) \\
\!\!\!\!\!\!\!\!\!\!\!\!\!\!\! && \quad
+ \; C(p,\rho_0, T) \Big( \frac 1N \sum_{x\in \bb Z}
\big\vert H(x/N) \big\vert \Big)^2 
\end{eqnarray*}
for some finite constant $C(p,\rho_0, T)$. Here $F(a) = a(1-a)$.

\section{ Gaussian tail estimates for labeled exclusion processes}
\label{sec3}

Fix $n\ge 2$ and a finite range, symmetric and irreducible transition
probability $p(\cdot)$ on $\bb Z^d$.  Consider $n$ labeled particles
moving on the $d$-dimensional lattice $\bb Z^d$ through stirring. This
dynamics can be informally described as follows. The $n$ particles
start from $n$ distinct sites of $\bb Z^d$. For each pair $(x,y)$ of
$\bb Z^d$, at rate $p(y-x)$, particles at $x$, $y$ exchange their
positions. This means that if there is a particle at $x$ (resp. $y$)
and no particle at $y$ (resp. $x$), the particle jumps from $x$ to $y$
(resp. from $y$ to $x$). If both sites are occupied, the particles
change their position and if none of them are occupied, nothing
happens.

The state space of this Markov process, denoted by $\mc B_n$, consists
of all vectors $\mb x=(x_1, \dots, x_n)$ of $(\bb Z^d)^n$ with
distinct coordinates:
$$
\mc B_n \;=\; \Big\{ \mb x=(x_1, \dots, x_n) \in (\bb Z^d)^n \,:\, 
x_i \not =x_j \text{ for } i\not = j\Big\}
$$
while the generator $L_n$ is given by
\begin{equation}
\label{eq:1}
(L_nf)(\mb x) \;=\; \sum_{x,y\in\bb Z^d} p(y-x) [f(\sigma^{x,y} \mb x) 
- f(\mb x)] \; .
\end{equation}
In this formula, for a configuration $\mb x = (x_1, \dots, x_n)$ in
$\mc B_n$,\,\, $\sigma^{x,y}\, \mb x$ is the configuration defined by
\begin{equation*}
(\sigma^{x,y}\, \mb x)_i \;=\; 
\left\{
\begin{array}{ll}
y & \text{if $x_i=x$, } \\
x & \text{if $x_i=y$, } \\
x_i & \text{otherwise.}
\end{array}
\right.
\end{equation*}
This generator corresponds to the generator \eqref{eq:9} in which
particles have been labeled and are therefore distinguishable.

It is easy to check that the counting measure on $\mc B_{n}$, denoted
by $\mu_{n}$, is an invariant reversible measure for the process.  The
goal of this section is to obtain sharp estimates on the transition
probability of this Markov process. To state the main results of the
section, fix a state $\mb z$ in $\mc B_n$ and denote by $f_t$ the
solution of the forward equation:
\begin{equation}
\label{eq:01}
\left\{
\begin{array}{l}
\partial_t f_t = L_n^* f_t \; , \\
f_0 (\mb x) = \mb 1\{\mb x = \mb z\} \; .
\end{array}
\right.
\end{equation}

Recall that we denote by $\Phi$ the Legendre transform of the convex
function $w^2 \cosh w$.

\begin{theorem} 
\label{s1}
Fix $n\ge 1$ and a point $\mb z = (z_1,\dots, z_n)$ in $\mc B_n$.  Let
$f_t$ be a solution of the forward equation (\ref{eq:01}). There
exist finite constants $C_2=C_2(n,d,p)$, $a_0=a_0(p)$ such that
$$
f_T(\mb x) \le \frac{C_2}{(1+T)^{nd/2}}\,\exp\Big\{
- \frac { a_0 T}{2 (\log T)^2} \, \Phi\Big( 
\frac{\Vert \mb x -\mb z \Vert \log T} {a_0^2  T}\Big)\Big\}
$$
for every $T > C_2$ and every configuration $\mb x$.
\end{theorem}

Since $\Phi (w) \sim w^2$ for $w$ small, for $\gamma>0$, there exists
a finite constant $a_1=a_1(p,\gamma)$ such that
$$
f_T(\mb x) \le \frac{C_2}{(1+T)^{nd/2}}\,\exp\Big\{\frac
{- \Vert \mb x - \mb z \Vert^2}{a_1 T} \Big\}
$$
for every $T > C_2$ and every configuration $\mb x$ such that
$\Vert \mb x - \mb z \Vert \le \gamma T / \log T$. 

On the other hand, since $x^2 \cosh x \le 2 e^{2x}$,
$\Phi(u) \ge (u/2) \log (u/4 e)$. Hence,
$$
f_T(\mb x) \le \frac{C_2}{(1+T)^{nd/2}}\,\exp\Big\{-\, \frac
{\Vert \mb x - \mb z \Vert}{ 4 a_0 \log T} \, \log \frac
{\Vert \mb x - \mb z \Vert \log T}{4 e a_0^2 T}\Big\}
$$
for every $T > C_2$. Of course this estimate is only interesting if
$\Vert \mb x - \mb z \Vert \gg  T / \log T$,

Since the proof of Theorem \ref{s1} follows closely the one of Theorem
2.2 in \cite{l}, we present only the main differences. Throughout this
section, $C_0$ stands for a universal constant, which may change from
line to line.

We first need a logarithmic Sobolev inequality for the process $X_t$
restricted to cubes. Fix an integer $\ell$ and decompose the lattice
$\bb Z^d$ into disjoint cubes $\{\Lambda_k \,:\, k \ge 1\}$ of length
$\ell$:
\begin{eqnarray*}
\!\!\!\!\!\!\!\!\!\!\!\!\! &&
\Lambda_k \;=\;  x_k + \{1,\dots,\ell\}^d \text{ for some $x_k$ in
  $\bb Z^d$} \,; \\
\!\!\!\!\!\!\!\!\!\!\!\!\! && \quad 
\Lambda_k \cap \Lambda_j \;=\; \phi 
\text{ for } k \ne j \quad
\text{and } \bigcup_{k\ge 1} \Lambda_k \;=\; \bb Z^d \,.
\end{eqnarray*}
For a vector $\mb k = (k_1, \dots, k_n)$, let $\Lambda_{\mb k}$ be the
finite cube of $(\bb Z^d)^n$ defined by $\Lambda_{\mb k} =
\Lambda_{k_1} \times \cdots \times \Lambda_{k_n}$ and let
$L_{\Lambda_{\mb k}}$ be the generator $L_n$ introduced in
(\ref{eq:1}) restricted to the cube $\Lambda_{\mb k}$. This means that
jumps from $\Lambda_{\mb k}$ to its complement are forbidden as well as
jumps from the complement to $\Lambda_{\mb k}$. 

\begin{lemma}
\label{s2}
There exists a finite constant $C_1$ depending only on the transition
probability $p(\cdot)$, the dimension $d$ and the total number of
particles $n$ such that
\begin{equation}
\label{eq:2}
\sum_{\mb x\in \Lambda_{\mb k}} f(\mb x) \log f(\mb x) 
\;\le\; C_1 \ell^2  \sum_{\mb x, \mb y \in \Lambda_{\mb k}} 
\Big\{ \sqrt{f(\mb y)} - \sqrt{f(\mb x)}\Big\}^2\; ,
\end{equation}
for all densities $f$ with respect to the uniform probability measure
over $\Lambda_{\mb k}$. In this formula, the sum on the right hand
side of the inequality is carried over all pairs $\mb x$, $\mb y$ in
$\Lambda_{\mb k}$ such that $\mb y = \sigma^{x,y} \mb x$ for some
$x$, $y$ with $p(y-x)>0$.
\end{lemma}

\begin{proof}
It is well known that a symmetric random walk evolving on a
$d$-dimensional cube satisfies a logarithmic Sobolev inequality of
type (\ref{eq:2}) and that the superposition of independent processes
satisfying logarithmic Sobolev inequalities also satisfies a
logarithmic Sobolev inequality, the constant being the maximum of the
individual constants. This proves (\ref{eq:2}) in the case where the
cubes $\Lambda_k$ are all different: $k_i \not = k_j$ for $i \not =
j$.

It remains to consider the case where some cubes are equal. In this
situation, the diagonal is forbidden because two particles cannot
occupy the same site, and two particles may exchange their position.
Fix $2\le m\le n$ and consider the hypercube $\Lambda_{\mb k} =
\Lambda_k \times \cdots \times \Lambda_k$ of $(\bb Z^d)^m$. If we do
not distinguish particles, we retrieve the symmetric simple exclusion
process on $\Lambda_k$ with $m$ particles.  This process satisfies a
logarithmic Sobolev inequality of type (\ref{eq:2}) with a constant
$C_0$ depending only on the dimension $d$ and the transition
probability $p(\cdot)$ \cite{y}. It is not difficult to recover
(\ref{eq:2}) for the random walk $X_t$ on $\Lambda_{\mb k}$ from this
estimate.

Indeed, let $\Sigma_{\Lambda_k,m}$ be the subsets of $\Lambda_k$ with
$m$ points: $\Sigma_{\Lambda_k,m} =\{A\subset \Lambda_k \,:\,
|A|=m\}$, let $\mu_{\Lambda_k, m}$ be the uniform probability measure
on $\Sigma_{\Lambda_k,m}$ and, for a density $f:\Lambda_{\mb k} \to\bb
R_+$ with respect to the uniform measure over $\Lambda_{\mb k}$, let
$\tilde f: \Sigma_{\Lambda_k,m} \to \bb R_+$ be the density with
respect to $\mu_{\Lambda_k, m}$ defined by
$$
\tilde f (\{x_1, \dots, x_m\})\;=\; \frac 1{m!} \sum_{\sigma} 
f(x_{\sigma(1)}, \dots, x_{\sigma (m)})\; ,
$$
where the summation is performed over all permutations $\sigma$ of
$m$ elements. 

With this notation, we may rewrite the left hand side of (\ref{eq:2}) as
\begin{equation}
\label{eq:3}
\sum_{A\in \Sigma_{\Lambda_k,m}} \tilde f (A) \sum_{\mb x \in A}
\frac {f(\mb x)}{\tilde f (A)} \log \frac {f(\mb x)}{\tilde f (A)}
\;+\; m! \sum_{A\in \Sigma_{\Lambda_k,m}} \tilde f (A) \,
\log {\tilde f (A)}\; ,
\end{equation}
where the summation over $\mb x$ is carried over all points $\mb x
=(x_1, \dots, x_m)$ such that $\{x_1, \dots, x_m\} = A$. 

It is not difficult to prove a logarithmic Sobolev inequality for the
permutation of $m$ points. Let $\bb S_m$ be the set of all
permutations $\sigma$ of $m$ points. Consider the Dirichlet form
$D_{\bb S_m}$ defined by
$$
D_{\bb S_m} (g) \;=\; \sum_{\sigma,\tilde \sigma\in \bb S_m}
[ g(\sigma) - g(\tilde \sigma)]^2\; .
$$
There exists a finite constant $C_0$ such that
$$
\sum_{\sigma\in \bb S_m} g(\sigma) \, \log g(\sigma)
\;\le\; C_0 D_{\bb S_m} (\sqrt{g})
$$
for all densities $g$ with respect to the uniform probability
measure on $\bb S_m$.

Since $f(\mb x) / \tilde f (A)$ is a density with respect to the
uniform probability measure over the set of all permutations, the
first term is bounded above by
\begin{equation}
\label{eq:4}
C_0 \sum_{A\in \Sigma_{\Lambda_k,m}} \sum_{\mb x, \mb y \in A}
\Big\{ \sqrt{f(\mb y)} - \sqrt{f(\mb x)} \Big\}^2
\end{equation}
for some finite universal constant. It remains to connect each point
$\mb x$ in $A$ to each point $\mb y$ in $A$ by a path $\mb x = \mb
z_0, \dots, \mb z_r= \mb y$ such that $\mb z_{j+1} = \sigma^{x,y} \mb
z_j$ for some $x$, $y$ with $p(y-x)>0$ to estimate the previous term
by the right hand side of (\ref{eq:2}). This can be done as follows.

Assume first that $d=1$. To explain the strategy in a simple way, we
allow two particles to occupy the same site in the construction of the
path. The modifications needed to respect the exclusion rule are
straightforward.  Fix $\mb x$ and $\mb y$ in a same set $A$.  Since
both points belong to the same set, there exists a permutation
$\sigma$ of $m$ points such that $y_i = x_{\sigma (i)}$ for $1\le i\le
m$. The path $\{\mb z_j\}$ connecting $\mb x$ to $\mb y$ is defined as
follows. We start changing the first coordinate $x_1$ of $\mb x$,
keeping all the other constants, moving from $\mb x =(x_1, \dots,
x_m)$ to $\mb w_1 = (y_1=x_{\sigma(1)}, x_2, \dots, x_m)$.  Note that
the last configuration has two particles occupying the same site. At
this point, we change the coordinate $x_{\sigma (1)}$, moving from a
new configuration $\mb w_2$, which is obtained from $\mb x$, by
replacing $x_1$ by $x_{\sigma(1)}$ and $x_{\sigma(1)}$ by
$x_{\sigma^2(1)}$, where $\sigma^2 = \sigma \circ \sigma$. We repeat
this procedure. If the orbit of $1$ for the permutation $\sigma$ is
the all set $\{1, \dots, m\}$, this algorithm produces a path from
$\mb x$ to $\mb y$. Otherwise, after completing the orbit of $1$ by
the map $\sigma$, we choose the smallest coordinate not belonging to
the orbit of $1$ and repeat the procedure.

Denote by $\Gamma_{\mb x, \mb y}$ the path just constructed.  Notice
that 
\renewcommand{\labelenumi}{\theenumi .}
\begin{enumerate}
\item its length is bounded by $m \ell$ and
  
\item all coordinates but one of each site $\mb z$ in $\Gamma_{\mb x,
    \mb y}$ belong to the set $\{x_1, \dots , x_m\}$.
\end{enumerate}

Therefore, by Schwarz inequality, (\ref{eq:4}) is bounded above by
\begin{eqnarray*}
&& C_0 \sum_{A\in \Sigma_{\Lambda_k,m}} \sum_{\mb x, \mb y \in A}
\big\vert \Gamma_{x,y} \big\vert \sum_{b\in\Gamma_{\mb x,\mb y}}
\Big\{ \sqrt{f(\mb b_2)} - \sqrt{f(\mb b_1)} \Big\}^2 \\
&& \quad \le\; C_0 \ell m \sum_{b} \Big\{ \sqrt{f(\mb b_2)} 
- \sqrt{f(\mb b_1)} \Big\}^2
\sum_{A\in \Sigma_{\Lambda_k,m}} \sum_{\substack{\mb x, \mb y \in A\\
  b\in \Gamma_{\mb x,\mb y}}} \;.
\end{eqnarray*}
The last sum in the first line is performed over all pairs $b
=(b_1,b_2)$ of consecutive sites in the path $\Gamma_{\mb x,\mb y}$,
while the first sum in the second line is performed over all pairs
$b=(b_1,b_2)$ such that $b_2 = \sigma^{x,x+y}$ for some $x$, $y$ in
$\bb Z^d$ such that $p(y)>0$. Since all but one coordinate of each
site in $\Gamma_{\mb x,\mb y}$ belong to $\{x_1, \dots, x_m\}$, for
each fixed bond $b=(b_1,b_2)$ there is at most $m \ell$ possible sets
$A$ which might use this bond.  For each set $A$, there is at most
$m!$ end points and $m!$ starting points for the path. The last sum is
thus bounded by
$$
\le\; C_0 \ell^2 m^2 (m!)^2  \sum_{b} \Big\{ \sqrt{f(\mb b_2)} 
- \sqrt{f(\mb b_1)} \Big\}^2\; .
$$
This concludes the proof of the estimate of the first term in
(\ref{eq:3}) in dimension 1.

The proof in higher dimension is similar. The idea is to consider a
configuration $\mb x$ as a point in $\bb Z^{md}$ and repeat the
previous algorithm, moving the first coordinate of the first particle,
then moving the first coordinate of the $\sigma (1)$-particle, until
all first coordinates of all particles are modified. At this point, we
change the second coordinate of the first particle and repeat the
procedure. This method gives a path of length at most $C_0 \ell m d$
and whose sites have all but one of the $md$ coordinates equal to the
coordinates of $\mb x$. These two properties permit to derive the
estimates obtained in dimension $1$, replacing $m$ by $md$. This
proves that the first term in (\ref{eq:3}) is bounded above by the the
right hand side of (\ref{eq:2}).

We focus now on the second term of (\ref{eq:3}).  By the logarithmic
Sobolev inequality for $m$ exclusion particles in a cube $\Lambda_k$,
this expression is less than or equal to
$$
C m! \sum_{A\in \Sigma_{\Lambda_k,m}} \sum_{x, y\in\bb Z^d} p(y) 
\Big\{ \sqrt{\tilde f (A_{x,x+y})} - \sqrt{\tilde f (A)} \Big\}^2
$$
for some finite constant $C$ depending only on $p(\cdot)$ and $d$. By
Schwarz inequality, this expression is bounded by the right hand side
of (\ref{eq:2}).
\end{proof}

The second main ingredient in the proof of Theorem \ref{s1} is an
estimate of the action of the generator $L_n$ on certain exponential
functions.

For a vector $\theta = (\theta_1, \dots, \theta_n)$, $\theta_i$ in
$\bb R^d$, denote by $\psi_\theta$ the function $\psi_\theta \colon
\mc B_n \to \bb R$ defined by $\psi_\theta (\mb x) = \exp\{\theta
\cdot \mb x \}$. Here, $\mb x \cdot \mb y$ represents the inner
product in $(\bb Z^d)^n$. An elementary computation shows that there
exists a finite constant $a_0$, depending only on the transition
probability $p(\cdot)$, such that
\begin{equation}
\label{eq:02}
\frac{1}{\psi_\theta (\mb x)} (L_n \psi_\theta )(\mb x) \; \le \; 
R(\theta)\;, \quad 
\sum_{x,y\in\bb Z^d} p(y-x) \Big\{ \frac {\psi_\theta (\sigma^{x,y} \mb x)}
{\psi_\theta (\mb x) } - 1 \Big\}^2 \;\le\; R(\theta) 
\end{equation}
for all $\mb x$ in $\mc B_n$, where
\begin{equation}
\label{eq:8}
R(\theta) \;=\; a_0 \sum_{j=1}^d\sum_{i=1}^n 
\Big( \cosh \{a_0 \theta_{i,j} \} - 1 \Big)\;.
\end{equation}

Next result relies mainly on Lemma \ref{s2} and on the bounds
(\ref{eq:02}). Its proof follows closely the one of Lemma 4.3 in
\cite{l} and is therefore omitted. For a positive function $\psi : \mc
B_n\to\bb R$, denote by $\mf D_{\psi}$ the Dirichlet formula defined
by
$$
\mf D_\psi (u) = (1/2) \sum_{\mb x \in \mc B_n} \sum_{x, y \in \bb
  Z^d} p(y-x) \big\{ u(\sigma^{x,y} \mb x)- u(\mb x) \big\}^2 \,
\psi(\mb x) \; .
$$

\begin{lemma}
\label{s3}
Fix a vector $\theta$ in $(\bb R^d)^n$, $\ell\ge 2$, denote by
$C_1$ the constant introduced in Lemma \ref{s2} and let $\psi =
\psi_\theta$.  There exists a finite constant $a_0$, depending only on
the transition probability $p(\cdot)$, such that
$$
\int f \, \log f\, \psi \,  d\mu_n 
\;\le\; - \int f \, \psi \, \log\psi \, d\mu_n 
\;-\; \log |\Lambda_\ell |^n \;+\;
4 C_1 \ell^2 \mf D_{\psi}(\sqrt{f}) \;+\; R(\theta)\ell^2
$$
for every density $f:\mc B_n\to\bb R$ with respect to $\psi \,
d\mu_n$.  
\end{lemma}

The estimates (\ref{eq:02}) permit also to prove the following bound.
Recall that $f_t$ is the solution of the forward equation
(\ref{eq:01}) and that $\mu_n$ is the counting measure on $\mc B_n$.

\begin{lemma}
\label{s4}
Fix a smooth increasing function $p: \bb R_+\to (1,\infty)$ and a
smooth function $\lambda = (\lambda_1, \dots, \lambda_n): \bb
R_+\to(\bb R^d)^n$.  Let $\psi_t (\mb x) = \exp\{ \lambda (t) \cdot
\mb x \}$ and let $h_t = f_t/\psi_t$, $u_t=h_t^{p(t)/2}$. There exists
a finite constant $a_0$, depending only on the transition probability
$p(\cdot)$, such that
\begin{eqnarray*}
\!\!\!\!\!\!\!\!\!\!\!\!\!\!
\frac{d}{dt} \int h_t^{p(t)} \, \psi_t\, d\mu_n & \le &
\frac{\dot p(t)}{p(t)} \int u_t^2 \, \log u_t^2 \, \psi_t \,d\mu_n -
(p(t)-1) \int u_t^2 \, \dot \psi_t \, d\mu_n \\
&-&\frac{(p(t)-1)}{p(t)} \mf D_{\psi_t}(u_t) \;+\; 
R(\lambda(t)) \, p(t) \int u_t^2\, \psi_t\,d\mu_n\; .
\end{eqnarray*}
\end{lemma}

The proof of Lemma \ref{s4} relies on the estimates (\ref{eq:02})
and follows closely the proof of Lemma 5.1 in \cite{l}.  
\medskip

We are now in a position to prove Theorem \ref{s1}.  Recall that $f_t$
is the solution of the forward equation (\ref{eq:01}).  Fix $T>0$
large, set $q= 1 + (\log T)^{-1}$, $q' = \log T$ and consider a smooth
increasing function $p\colon [0,T]\to [q,q']$ such that $p(0)=q$,
$p(T)=q'$. At the end of the proof, $p(t)$ will be taken as a
rescaling of the function $[1-(s/T)^\alpha]^{-1}$ for some
$0<\alpha<1/2$.

Following Davies \cite{d}, fix $\theta = (\theta_1, \dots,
\theta_n)$ in $(\bb R^d)^n$, define $\psi_t\colon \mc B_n \to \bb R_+$
by
$$
\psi_t(\mb x) \;=\; \exp\Big\{\frac{p(t)}{p(t)-1} \,
\theta \cdot \mb x \Big\}\; ,
$$
denote $\theta_i\,p(t)/[p(t)-1]$ by $\lambda_i(t)$ and let $h_t =
f_t/\psi_t$.  

For a function $g\colon \mc B_n \to \bb R$ and $1\le
p<\infty$, denote by $\Vert g \Vert_{\psi, p}$ the $L^p$ norm of $g$
with respect to the measure $\psi \, d\mu_n$~:
$$
\Vert g \Vert_{\psi, p}^p \; =\; \sum_{\mb x \in\mc B_n}
\vert g(\mb x)\vert ^p \psi (\mb x) \; .
$$
A straightforward computation gives that
\begin{equation}
\label{eq:5}
\frac{d}{dt} \log \Vert h_t \Vert_{\psi_t, p(t)} \;=\; 
-\frac{\dot p(t)}{p(t)} \log \Vert h_t\Vert_{\psi_t, p(t)} 
\;+\; \frac{1}{p(t)} \frac{1} {\Vert h_t\Vert_{\psi_t, p(t)}^{p(t)}} 
\frac{d}{dt} \Vert h_t\Vert_{\psi_t, p(t)}^{p(t)}\;.
\end{equation}
Denote $h_t^{p(t)}$ by $u_t^2$ and $u_t^2/\Vert u_t\Vert_2^2$ by
$v_t^2$. By Lemma \ref{s4}, the second term on the right hand side of
last formula is bounded above by
\begin{eqnarray}
\label{eq:6}
\!\!\!\!\!\!\!\!\!\!\!\!\!\! && 
\frac{\dot p(t)}{p(t)^2} \int v_t^2 \, \log v_t^2 \, \psi_t\,d\mu_n 
\;+\; \frac{\dot p(t)}{p(t)} \log\Vert h_t\Vert_{\psi_t, p(t)} 
\;+\; R(t) \\
\!\!\!\!\!\!\!\!\!\!\!\!\!\! && \quad 
- \;\frac{p(t)-1}{p(t)^2}\, \mf D_{\psi_t}(v_t) \;-\; 
\frac{p(t)-1}{p(t)} \int v_t^2 \, \dot \psi_t\, d\mu_n\; ,
\nonumber
\end{eqnarray}
where $R(t)=R(\lambda(t))$.  Notice that the second term in this
expression cancels with the first term in the previous formula and
that $v_t^2$ is a density with respect to the measure
$\psi_t\,d\mu_n$\,. By Lemma \ref{s3}, the first term of this formula
is bounded by
\begin{equation}
\label{eq:7}
\frac{\dot p(t)}{p(t)^2} \Big\{ - \int v_t^2 \, \psi_t\, \log\psi_t\, 
d\mu_n \;-\; \log |\Lambda_\ell |^n \;+\;  4 C_1 \ell^2 \mf D_{\psi_t}(v_t)
\;+\; \ell^2 R(t) \Big\}
\end{equation}
for all $\ell \ge 2$.  

By definition of $\psi_t$,
$$
-\frac{p(t)-1}{p(t)}\, \dot \psi_t \;=\;
\frac{\dot p(t)}{p(t)^2}\, \psi_t \log \psi_t\,,
$$
so that the first term of formula (\ref{eq:7}) cancels with the
fifth term of formula (\ref{eq:6}). Denote by $[a]$ the integer part
of a real $a$. If we set $\ell = \ell(t)$ as
$$
\ell(t) = \Big[ \,\sqrt{\frac{p(t)-1}{4 C_1 \dot p(t)} }\,
\Big] \, , 
$$
a straightforward computation shows that the Dirichlet form in
formula (\ref{eq:6}) cancels with the Dirichlet form appearing in
(\ref{eq:7}). The inequality $\ell (t)\ge 2$ imposes conditions on
$p(t)$ that will need to be checked when defining $p(t)$.

Up to this point we proved that
$$
\frac {d}{dt} \log \Vert h_t\Vert_{\psi_t,p(t)} \;\le\;
- \frac{\dot p(t)}{p(t)^2} \log |\Lambda_\ell |^n \; +\; 2 R(t)
$$
because $\ell(t)^2 \dot p(t) \le p(t)^2$ by definition of $\ell
(t)$. Integrating in time, we obtain that
\begin{equation*}
\Vert h_T \Vert_{\psi_T, p_T} \;\le\; \Vert h_0 \Vert_{\psi_0, p_0}
\exp\Big\{ -(nd/2) \int_0^T dt\, \frac {\dot p(t)}{p(t)^2} \log 
\frac {p(t) -1}{8 C_1 \dot p(t)} \;+\; 2 \int_0^T dt\, R(t) \Big\}\; .
\end{equation*}
because $\ell(t)^2 \ge [p(t)-1]/8 C_1 \dot p(t)$.  By definition of
the density $f$, $\Vert h_0 \Vert_{\psi_0, p_0} = f(\mb z) \psi_0(\mb
z)^{1-p_0/p_0}= \exp \{\theta \cdot \mb z\}$. On the other hand,
$\Vert h_T \Vert_{\psi_T, p_T}$ is bounded below by $f_T(\mb x)
\psi_T(\mb x)^{1-p_T/p_T} = \exp\{\theta \cdot \mb x\}$ for every $\mb
x$ in $\mc B_n$. Moreover, since $- \dot p(t)/p(t)^2 = (1/p(t))'$,
$$
- \log (8 C_1)  \int_0^T dt\, \frac {\dot p(t)}{p(t)^2}
\;\le\; \log (8 C_1) 
$$
because $p(0) = 1 + (\log T)^{-1}$, $p(1) = \log T$. Finally, since
$p(t)$ is an increasing function, $p(t)/[p(t)-1] \le 1 + \log T\le 2
\log T$ for $T\ge e$.  Therefore, if we assume that $\vert\theta_{i,j}
\vert \le B/2\log T$ for some finite constant $B$, $R(t) \le C(a_0, B)
\Vert\theta\Vert^2 p(t)^2/(p(t)-1)^2$, where $C(a_0, B) = 2 a_0^3
M(a_0 B)$ and
$$
M(r) \;:=\; \sup_{|w|\le r} \frac{\cosh w -1}{w^2} \;\le\; \cosh r\; .
$$

Putting together all previous estimates, we obtain that
\begin{eqnarray*}
\!\!\!\!\!\!\!\!\!\!\!\!\!\!\! &&
f_T(\mb x) \;\le\;  C(n,d) \exp\{ - \theta
\cdot (\mb x - \mb z)\} \, 
\exp\Big\{ -(n d/2) \int_0^T dt\, \frac {\dot p(t)}{p(t)^2} \log 
\frac {p(t) -1}{\dot p(t)} \Big\} \;\times \\
\!\!\!\!\!\!\!\!\!\!\!\!\!\!\! &&
\qquad\qquad\qquad\qquad\qquad \times\; 
\exp\Big\{C(a_0, B) \Vert\theta\Vert^2 
\int_0^T dt\, \frac {p(t)^2} {[p(t) -1]^2} \Big\}
\end{eqnarray*}
provided $\vert\theta_{i,j} \vert \le B/2\log T$.

It remains to choose an appropriate increasing smooth function
$p:[0,T] \to [q,q']$ which connects $1 + (\log T)^{-1}$ to $\log T$ to
conclude the proof of the theorem. Let $q(s)= p(sT)/p(sT)-1$ and
notice that $q(0) = \log T + 1$, $q(1)=\log T/\log T - 1$. With this
notation, a change of variables and an elementary computation shows
that the two previous integrals become
\begin{eqnarray*}
\!\!\!\!\!\!\!\!\!\!\!\!\!\!\! &&
(nd/2) \int_0^1 ds\, \frac{q'(s)}{q(s)^2} \log
\frac{q(s)-1}{-q'(s)} \;-\; 
(nd/2) \log T \Big( \frac{\log T}{\log T+1} - \frac 1{\log T}\Big) \\
\!\!\!\!\!\!\!\!\!\!\!\!\!\!\! && \quad 
\text{and}\quad 
C(a_0, B) \Vert\theta\Vert^2 T \int_0^1 ds\, q(s)^2 \; .
\end{eqnarray*}
The second term of the first line is bounded by $\log T^{-nd/2} +
C(n,d)$, which is responsible for the diagonal estimate of the
density. 

Let $g(s) = s^{-\alpha}$ for some $0<\alpha<1/2$. It is easy to show 
\begin{equation*}
\int_0^1 ds\, \frac{- g'(s)}{g(s)^2} \log \frac{g(s)-1}{-g'(s)} \;<\; 
\infty \;, \quad  \int_0^1 ds\, g(s)^2 \;<\; \infty \; .
\end{equation*}
Defining $q(s) = g(a +(b-a)s)$ for appropriate constants $a$, $b$, we
deduce that
$$
f_T(\mb x) \;\le\;  \frac{C(n,d)}{T^{nd/2}} 
\exp\{ - \theta \cdot (\mb x - \mb z) + C(a_0, B) \Vert\theta\Vert^2 T\} 
$$
provided $\vert\theta_{i,j} \vert \le B/2\log T$. An elementary
computation shows that with this choice $\ell(t)\ge 2$ for all $0\le
t\le T$ provided $T$ is chosen large enough: $T\ge C_2(n,d)$.

Fix $\mb x$, let $\mb y = \mb x - \mb z$ and choose $\theta = B
(2 \log T)^{-1} \mb y/\Vert \mb y\Vert$ so that $|\theta_{i,j}|\le
B/2 \log T$. With this choice, the expression inside braces in the
previous formula becomes bounded by
$$
- \frac {\Vert \mb x - \mb z \Vert \, B}{2 \log T} \;+\; 
\frac{ C(a_0,B) B^2 T} {(2\log T)^2}
$$
for every $B>0$. Recall the definition of $C(a_0, B)$, change
variables as $B' = a_0 B$ and minimize over $B'$ to obtain that the
previous expression is bounded above by
$$
- \frac { a_0 T}{2 (\log T)^2} \, \Phi\Big( \frac{\Vert \mb x -\mb z \Vert 
\log T} {a_0^2  T}\Big)\;,
$$
where $\Phi$ is the convex conjugate of $w^2 \cosh w$. This
concludes the proof of Theorem \ref{s1}.

\medskip
\noindent{\bf Proof of Theorem \ref{s5}}.
Theorem \ref{s5} follows from Theorem \ref{s1} since the evolution of
$n$ random walks evolving with exclusion can be obtained from the
evolution of $n$ labeled random walks by just ignoring the labels.

\end{document}